# A note on the estimation of extreme value distributions using maximum product of spacings

T. S. T. Wong[1] and W. K. Li[1]

*The University of Hong Kong*

**Abstract:** The maximum product of spacings (MPS) is employed in the estimation of the Generalized Extreme Value Distribution (GEV) and the Generalized Pareto Distribution (GPD). Efficient estimators are obtained by the MPS for all $\gamma$. This outperforms the maximum likelihood method which is only valid for $\gamma < 1$. It is then shown that the MPS gives estimators closer to the true parameters compared to the maximum likelihood estimates (MLE) in a simulation study. In cases where sample sizes are small, the MPS performs stably while the MLE does not. The performance of MPS estimators is also more stable than those of the probability-weighted moment (PWM) estimators. Finally, as a by-product of the MPS, a goodness of fit statistic, Moran's statistic, is available for the extreme value distributions. Empirical significance levels of Moran's statistic calculated are found to be satisfactory with the desired level.

## 1. Introduction

The GEV and the GPD (Pickands, [13]) distributions are widely-adopted in extreme value analysis. As is well known the maximum likelihood estimates (MLE) may fail to converge owing to the existence of an unbounded likelihood function. In some cases, MLE can be obtained but converges at a slower rate when compared to that of the classical MLE under regular conditions.

Recent studies (e.g. Juarez & Schucany, [9]) show that maximum likelihood estimation and other common estimation techniques lack robustness. In addition, the influence curve of the MLE is shown unstable when the sample size is small. Although new methods (Juarez & Schucany, [9]; Peng and Welsh, [12]; Dupuis, [6]) were proposed, arbitrary parameters are sometimes involved, resulting in more intensive computation which is in general undesirable. There have been studies in overcoming the difficulties of the MLE in extreme value analysis but none has considered the MPS. Furthermore, a goodness-of-fit test on the fitted GEV or GPD is rarely considered.

In this study, the MPS method will first be considered for the purpose of finding estimators which may not be obtained by the maximum likelihood method. As a by product, the Moran's statistic, a function of product of spacings, can be treated as a test statistics for model checking. This is one of the nice outcomes of MPS which Cheng and Stephens [4] demonstrated but is overlooked by the extreme value analysis literature.

In Section 2, we discuss some problems of the MLE. In Section 3, we formulate the MLE, the MPS and the Moran statistics. In Section 4, results of simulation studies are presented to evaluate the performance of the method proposed. In Section 5, we provide some real examples in which the MPS is more convincing. A brief discussion is presented in Section 6.

[1]Department of Statistics and Actuarial Science, The University of Hong Kong, Pokfulam Road, Hong Kong, e-mail: h0127272@hkusua.hku.hk

## 2. Problems of the MLE

The problems of the MLE in model fitting were discussed by Weiss and Wolfowitz [15]. Related discussions in connection to the Weibull and the Gamma distributions can be found in [2, 3, 5, 14]. Smith [14] found densities in the form

$$(2.1) \qquad f(x; \theta, \phi) = (x - \theta)^{\alpha - 1} g(x - \theta; \phi), \quad \theta < x < \infty$$

where $\theta$ and $\phi$ are unknown parameters and $g$ converges to a constant as $x \downarrow \theta$. As is well-known for $\alpha > 2$, the MLE is as efficient as in regular cases. For $\alpha = 2$, the estimated parameters are still asymptotically normal, but the convergence rate is $(n \log n)^{\frac{1}{2}}$ which is larger than the classical rate of $n^{\frac{1}{2}}$. For $1 < \alpha < 2$, the MLE exists but the asymptotic efficiency problem is not solved. And the order of convergence could be as high as $O(n^{\frac{1}{\alpha}})$. For $\alpha < 1$, MLE does not exist. Both the GEV and the GPD encounter the above difficulties as both can be reparameterised into the form (2.1).

As an alternative to the MLE, the MPS was established by Cheng and Amin [2]. With the MPS, not only can problems with non-regular condition be better solved, but models originally estimable under the MLE framework can also be better estimated by the MPS using a much simpler algorithm. Cheng & Amin [2] showed that the MPS estimators are asymptotically normal even for $0 < \alpha < 1$. This overcomes to a certain extent the weakness existing in the MLE. Hence, the MPS may be one of the most robust estimation techniques and yet the least computational expensive in extreme value analysis. The present paper employs the MPS in the estimation of the GEV and the GPD. On the other hand, many previous studies (Hosking, [7]; Marohn, [10]) concentrated on testing the shape parameter. Goodness-of-fit test on the model as a whole has been very few. In this study, the Moran's statistic (Cheng and Stephens, [4]; Moran, [11]) arising naturally as a by product of the MPS estimator was utilized to check the adequacy of the overall model.

## 3. Formulations of the MLE, the MPS and the Moran's statistic

### 3.1. The MLE and the MPS

The c.d.f of the GEV and the GPD are respectively

$$H(x; \gamma, \mu, \sigma) = \exp\left[-\left(1 - \gamma\frac{x-\mu}{\sigma}\right)^{\frac{1}{\gamma}}\right], \quad 1 - \gamma\frac{x-\mu}{\sigma} > 0;$$

and

$$G(x; \gamma, \sigma) = 1 - \left(1 - \gamma\frac{x}{\sigma}\right)^{\frac{1}{\gamma}}, \quad 1 - \gamma\frac{x}{\sigma} > 0.$$

where

$$\gamma \neq 0, \quad -\infty < \mu < \infty, \quad \sigma > 0$$

Let $h(x)$ and $g(x)$ be the corresponding densities,

$$h(x) = \frac{1}{\sigma}\left(1 - \gamma\frac{x-\mu}{\sigma}\right)^{\frac{1}{\gamma}-1} \exp\left[-\left(1-\gamma\frac{x-\mu}{\sigma}\right)^{\frac{1}{\gamma}}\right] ;$$

and

$$g(x) = \frac{1}{\sigma}\left(1 - \gamma\frac{x}{\sigma}\right)^{\frac{1}{\gamma}-1} .$$

The log-likelihood functions per observation are respectively

$$L_{\text{GEV}}(\gamma,\mu,\sigma) = -\log\sigma + \left(\frac{1}{\gamma}-1\right)\log\left(1-\gamma\frac{x-\mu}{\sigma}\right) - \left(1-\gamma\frac{x-\mu}{\sigma}\right)^{\frac{1}{\gamma}} ;$$

and

$$L_{\text{GPD}}(\gamma,\mu,\sigma) = -\log\sigma + \left(\frac{1}{\gamma}-1\right)\log\left(1-\gamma\frac{x}{\sigma}\right) .$$

Applying the same argument stated in [14], as $x \downarrow \mu + \frac{\sigma}{\gamma}$, the information matrix of $L_{\text{GEV}}(\gamma,\mu,\sigma)$ is infinite for $\gamma > \frac{1}{2}$. The same difficulty arises in the GPD as $x \downarrow \frac{\sigma}{\gamma}$. In this case, the underlying distribution is J-shaped where maximum likelihood is bound to fail. Worse still, MLEs (Denoted by $\hat{\boldsymbol{\Theta}}_{\text{GEV}} = (\hat{\gamma},\hat{\mu},\hat{\sigma})^T$ and $\hat{\boldsymbol{\Theta}}_{\text{GPD}} = (\hat{\gamma},\hat{\sigma})^T$ respectively for the GEV and the GPD) may not exist when $\gamma > 1$. Let $x_1 < x_2 < \cdots < x_n$ be an ordered sample of size $n$ and define spacings $D_i(\boldsymbol{\theta})$ by

$$\text{GEV}: D_i(\boldsymbol{\theta}) = H(x_i,\gamma,\mu,\sigma) - H(x_{i-1};\gamma,\mu,\sigma), \quad (i=1,2,\ldots,n+1);$$
$$\text{GPD}: D_i(\boldsymbol{\theta}) = G(x_i,\gamma,\sigma) - G(x_{i-1};\gamma,\sigma), \quad (i=1,2,\ldots,n+1);$$

where $H(x_0;\gamma,\mu,\sigma) \equiv G(x_0;\gamma,\sigma) \equiv 0$ and $H(x_{n+1};\gamma,\mu,\sigma) \equiv G(x_{n+1};\gamma,\sigma) \equiv 1$.

MPS estimators (Denoted by $\check{\boldsymbol{\Theta}}_{\text{GEV}} = (\check{\gamma},\check{\mu},\check{\sigma})^T$ and $\check{\boldsymbol{\Theta}}_{\text{GPD}} = (\check{\gamma},\check{\sigma})^T$ respectively for the GEV and the GPD) are found by minimizing

$$M(\boldsymbol{\theta}) = -\sum_{i=1}^{n+1} \log D_i(\boldsymbol{\theta}).$$

By taking the cumulative density in the estimation, the objective function $M(\boldsymbol{\theta})$ does not collapse for $\gamma < 1$ as $x \downarrow \mu + \frac{\sigma}{\gamma}$ for the GEV or as $x \downarrow \frac{\sigma}{\gamma}$ for the GPD. The MLE, however, does not have such an advantage. There is in probability a solution $\hat{\boldsymbol{\Theta}}$ that is asymptotically normal only for $\gamma < \frac{1}{2}$. The strength of MPS over MLE is demonstrated by the following two theorems.

**Theorem 3.1.** *Let $\boldsymbol{\Theta}_{0\text{GEV}} = (\gamma_0,\mu_0,\sigma_0)^T$ and $\boldsymbol{\Theta}_{0\text{GPD}} = (\gamma_0,\sigma_0)^T$ be the true parameters of the GEV and the GPD respectively. Under regularity conditions (See for example: [14])*

(i) *For $\gamma < \frac{1}{2}$, $n^{\frac{1}{2}}(\hat{\boldsymbol{\Theta}} - \boldsymbol{\Theta}_0) \xrightarrow{D} N\left(\boldsymbol{0}, -E\left(\frac{\partial^2 L}{\partial \boldsymbol{\Theta}^2}\right)^{-1}\right);$*

(ii) *For $\gamma = \frac{1}{2}$, $\left(\hat{\mu} + \frac{\hat{\sigma}}{\hat{\gamma}}\right) - (\mu_0 + \frac{\sigma_0}{\gamma_0}) \xrightarrow{D} O_p[(n\log n)^{-\frac{1}{2}}]$, and $n^{\frac{1}{2}}(\hat{\boldsymbol{\Theta}} - \boldsymbol{\Theta}_0) \xrightarrow{D} N(\boldsymbol{0}, -E(\frac{\partial^2 L}{\partial \boldsymbol{\Theta}^2})^{-1})$, where $\boldsymbol{\Theta} = (\gamma,\sigma)^T$;*

(iii) *For $\frac{1}{2} < \gamma < 1$, $\left(\hat{\mu} + \frac{\hat{\sigma}}{\hat{\gamma}}\right) - \left(\mu_0 + \frac{\sigma_0}{\gamma_0}\right) \xrightarrow{D} O_p(n^{-\gamma})$, and $n^{\frac{1}{2}}(\hat{\boldsymbol{\Theta}} - \boldsymbol{\Theta}_0) \xrightarrow{D} N\left(\boldsymbol{0}, -E\left(\frac{\partial^2 L}{\partial \boldsymbol{\Theta}^2}\right)^{-1}\right)$, where $\boldsymbol{\Theta}$ is as in (ii).*

(iv) *For $\gamma \geq 1$, the MLE does not exist.*

**Theorem 3.2.** *Under the same conditions as in Theorem 3.1*

(i) *For $\gamma < \dfrac{1}{2}$, $n^{\frac{1}{2}}(\breve{\Theta} - \Theta_0) \xrightarrow{D} N(\mathbf{0}, -E(\dfrac{\partial^2 L}{\partial \Theta^2})^{-1})$;*

(ii) *For $\gamma = \dfrac{1}{2}$, $(\breve{\mu} + \dfrac{\breve{\sigma}}{\breve{\gamma}}) - (\mu_0 + \dfrac{\sigma_0}{\gamma_0}) \xrightarrow{D} O_p[(n \log n)^{-\frac{1}{2}}]$, and $n^{\frac{1}{2}}(\breve{\Theta} - \Theta_0) \xrightarrow{D}$ $N(\mathbf{0}, -E(\dfrac{\partial^2 L}{\partial \Theta^2})^{-1})$, where $\Theta = (\gamma, \sigma)^T$;*

(iii) *For $\gamma > \dfrac{1}{2}$, $(\breve{\mu} + \dfrac{\breve{\sigma}}{\breve{\gamma}}) - (\mu_0 + \dfrac{\sigma_0}{\gamma_0}) \xrightarrow{D} O_p(n^{-\gamma})$, and $n^{\frac{1}{2}}(\breve{\Theta} - \Theta_0) \xrightarrow{D} N$ $(\mathbf{0}, -E(\dfrac{\partial^2 L}{\partial \Theta^2})^{-1})$, where $\Theta$ is as in (ii).*

Proofs of Theorems 3.1 and 3.2 follow the arguments in [14] and [2] respectively by checking the conditions therein.

It is obvious that efficient estimators can still be obtained by the MPS for $\gamma > \frac{1}{2}$ but not the MLE. From (iii) above, it is clear that the MPS still works while the MLE fails for $\gamma \geq 1$. It seems that it is a fact overlooked by researchers working in the extreme value literature.

### 3.2. Moran's statistic

In the MPS estimation, $M(\boldsymbol{\theta})$ is called the Moran's statistic which can be used as a test for a goodness-of-fit test. Cheng and Stephens [4] showed that under the null hypothesis, $M(\boldsymbol{\theta})$, being independent of the unknown parameters, has a normal distribution and a chi-square approximation exists for small samples with mean and variance approximated respectively by

$$\mu_M \approx (n+1)\log(n+1) - \frac{1}{2} - \frac{1}{12(n+1)},$$

and

$$\sigma_M^2 \approx (n+1)\left(\frac{\pi^2}{6} - 1\right) - \frac{1}{2} - \frac{1}{6(n+1)}.$$

Define

$$C_1 = \mu_M - \left(\frac{1}{2}n\right)^{\frac{1}{2}} \sigma_M, \quad C_2 = (2n)^{-\frac{1}{2}} \sigma_M.$$

The test statistic is

$$T(\breve{\boldsymbol{\theta}}) = \frac{M(\breve{\boldsymbol{\theta}}) + \frac{1}{2}k - C_1}{C_2}$$

which follows approximately a chi-square distribution of $n$ degrees of freedom under the null hypothesis. Monte Carlo simulation of the Weibull, the Gamma and the Normal distributions in [4] showed the accuracy of the test based on $T(\breve{\boldsymbol{\theta}})$. In the next section, we provide further evidence supporting the use of MPS for fitting the extreme value distributions.

### 4. Simulation study

A set of simulations was performed to evaluate the advantage of the MPS over the MLE of the GEV and the GPD based on selected parameters for different sample sizes $n = (10, 20, 50)$. Empirical significance levels of Moran's statistic were then considered using $\chi^2_{n,\alpha}$ as the benchmark critical value. Finally, data were generated from an exponential distribution and the cluster maxima of every 30 observations were fitted to the GEV.

The subroutine `DNCONF` in the IMSL library was used to minimize a function. The data analysed in the paper and the Fortran90 programs used in the computation are available upon request.

We have done extensive simulations to assess the performance of MPS estimators. Only four simulation results in each combination of $\gamma$ and $n$ are reported. The location and scale parameter, $\mu = 1$ and $\sigma = 1$, were used throughout. On the basis of the results from asymptotic normality of the MPS that were presented in Section 3, we chose a combination of $\gamma = (-0.2, 0.2, 1, 1.2)$ to compare the estimation performance between the maximum likelihood method and the MPS where the last two cases should break down for the MLE. 10000 simulations of sample sizes $n = (10, 20, 50)$ were performed. Data were generated from the same random seed and estimations were performed under the same algorithm. Define the mean absolute error for the MLE and the MPS respectively by

$$\frac{1}{l}|\hat{\boldsymbol{\Theta}}_l^T - \boldsymbol{\Theta}_0 \mathbf{1}^T|\mathbf{1} \quad \text{and} \quad \frac{1}{l}|\check{\boldsymbol{\Theta}}_l^T - \boldsymbol{\Theta}_0 \mathbf{1}^T|\mathbf{1} .$$

where $\hat{\boldsymbol{\Theta}}_l$ and $\check{\boldsymbol{\Theta}}_l$ are $l \times 1$ vectors of the MLE and MPS estimators respectively, $|\boldsymbol{Y}|$ means the element-wise absolute value of $\boldsymbol{Y}$, $p$ is the number of estimated parameters and $l = 10000$ is the number of replications. The mean absolute error measures the average deviation of estimators from the true parameters and hence is a measure of robustness. A small mean absolute error is expected.

As suggested by a referee, the MPS was also compared to the method of probability-weighted moment (PWM)(Hosking et al., [8]) for the GEV model. We followed Hosking's approach in his Table 3 and estimated the tail parameter by Newton-Raphson's Method. Tables 1 and 2 display the medians of the parameters in 10000 estimations together with the mean absolute error in bracket. Both the MPS estimates and the MLEs are in line with the true parameters but MPS tends to give a closer result for the GEV. It can also be seen that the MPS gives much more stable estimates than the MLE in general. For $\gamma = -0.2$ and $\gamma = 0.2$, the PWM performed well with slightly smaller mean absolute errors than the MPS. However, for $\gamma = 1$ and $\gamma = 1.2$, the bias of the PWM is rather severe. Note that some of the mean absolute errors for the MLE are unacceptably large due to serious outliers of estimated parameters. Non-regularity of the likelihood function caused occasional non-convergence. The frequency of such problems is reported in Tables 3 and 4. Failures of convergence were detected when the magnitudes of any estimator in an entry exceeds 100. The failure rates of MLE are relatively higher than those of MPS. Some estimated parameters of the MLE went up to as high as 500000. This explains the extremely large mean absolute errors of the MLE. Although there were failures in MPS, the maximum values were less than 1000, comparably less severe than the MLE. The PWM has zero failure rates but as mentioned above, it has a severe bias when $\gamma \geq 1$. It is noticed that the MLEs have smaller mean absolute error only in cases where sample size is large. However, the MPS estimators have virtually no fall off in its performance across sample sizes. These are in agreement with the theoretical results in Theorems 3.1 and 3.2. Overall, the MPS seems to be the most stable in its performance.

The Moran's statistic, $M(\boldsymbol{\theta})$, has a chi-square distribution with $n$ degrees of freedom. Monte Carlo simulations with 10000 observations per entry, each entry with sample size $n = (10, 20, 50)$ were conducted to compute the empirical significant levels. Again the null distributions were the models under consideration in Tables 1 and 2. It can be seen from Tables 5 and 6 that the empirical sizes for both the GEV and the GPD are very conservative at small sample sizes $n = 10$. Improvement was seen at $n = 20$. Though slightly conservative, it is acceptable in some applications. But the results at $n = 50$ are very good even with $\gamma = 1$ and $\gamma = 1.2$.

TABLE 1
*Simulation results of MPS estimates, MLEs and PWM estimates on the GEV. Shown are medians of estimated parameters from* 10000 *simulations of sample sizes* $n = (10, 20, 50)$. *Numbers in the bracket are mean absolute errors of estimates*

| $n$ | True parameters | | | MPS estimates | | | MLEs | | | PWM estimates | | |
|---|---|---|---|---|---|---|---|---|---|---|---|---|
| | $\gamma_0$ | $\mu_0$ | $\sigma_0$ | $\check{\gamma}$ | $\check{\mu}$ | $\check{\sigma}$ | $\hat{\gamma}$ | $\hat{\mu}$ | $\hat{\sigma}$ | $\gamma$ | $\mu$ | $\sigma$ |
| 10 | −0.2 | 1 | 1 | −0.28 | 0.96 | 1.06 | −0.22 | 1.09 | 0.98 | −0.15 | 1.02 | 0.93 |
| | | | | (0.38) | (0.30) | (0.30) | (435) | (175) | (295) | (0.18) | (0.29) | (0.25) |
| | 0.2 | 1 | 1 | 0.20 | 0.97 | 1.09 | 0.43 | 0.98 | 0.84 | 0.10 | 0.97 | 0.88 |
| | | | | (0.33) | (0.30) | (0.27) | (117) | (69) | (104) | (0.18) | (0.28) | (0.20) |
| | 1 | 1 | 1 | 1.17 | 0.98 | 1.11 | 1.20 | 0.78 | 0.78 | 0.59 | 0.89 | 0.88 |
| | | | | (0.53) | (0.33) | (0.45) | (54) | (90) | (50) | (0.42) | (0.31) | (0.27) |
| | 1.2 | 1 | 1 | 1.40 | 0.99 | 1.13 | 1.36 | 0.77 | 0.80 | 0.70 | 0.87 | 0.90 |
| | | | | (0.90) | (0.35) | (0.49) | (32) | (50) | (0.26) | (0.50) | (0.33) | (0.31) |
| 20 | −0.2 | 1 | 1 | −0.25 | 0.97 | 1.06 | −0.26 | 1.06 | 1.02 | −0.17 | 1.01 | 0.96 |
| | | | | (0.20) | (0.21) | (0.20) | (0.69) | (0.98) | (0.47) | (0.14) | (0.21) | (0.17) |
| | 0.2 | 1 | 1 | 0.20 | 0.98 | 1.07 | 0.25 | 1.09 | 1.00 | 0.14 | 0.98 | 0.94 |
| | | | | (0.17) | (0.20) | (0.17) | (46) | (1.01) | (50) | (0.13) | (0.20) | (0.14) |
| | 1 | 1 | 1 | 1.10 | 0.99 | 1.09 | 1.18 | 0.80 | 0.80 | 0.77 | 0.94 | 0.95 |
| | | | | (0.36) | (0.24) | (0.27) | (85) | (34) | (95) | (0.27) | (0.21) | (0.20) |
| | 1.2 | 1 | 1 | 1.33 | 0.99 | 1.10 | 1.35 | 0.79 | 0.81 | 0.92 | 0.93 | 0.96 |
| | | | | (0.55) | (0.26) | (0.31) | (13) | (50) | (0.33) | (0.32) | (0.22) | (0.22) |
| 50 | −0.2 | 1 | 1 | −0.22 | 0.99 | 1.04 | −0.25 | 1.02 | 1.0 | −0.18 | 1.01 | 0.98 |
| | | | | (0.11) | (0.13) | (0.11) | (0.12) | (0.13) | (0.11) | (0.10) | (0.13) | (0.11) |
| | 0.2 | 1 | 1 | 0.20 | 0.99 | 1.04 | 0.20 | 1.04 | 1.01 | 0.18 | 0.99 | 0.97 |
| | | | | (0.09) | (0.13) | (0.10) | (0.09) | (0.13) | (0.09) | (0.08) | (0.13) | (0.09) |
| | 1 | 1 | 1 | 1.05 | 0.99 | 1.05 | 1.11 | 0.89 | 0.88 | 0.90 | 0.97 | 0.98 |
| | | | | (0.29) | (0.16) | (0.16) | (50) | (50) | (0.22) | (0.16) | (0.13) | (0.12) |
| | 1.2 | 1 | 1 | 1.27 | 0.99 | 1.06 | 1.29 | 0.87 | 0.89 | 1.08 | 0.97 | 0.99 |
| | | | | (0.14) | (0.12) | (0.17) | (0.15) | (0.15) | (0.15) | (0.19) | (0.13) | (0.14) |

TABLE 2
*Simulation results of MPS estimates and MLEs on the GPD. Shown are medians of estimated parameters from* 10000 *simulations of sample sizes* $n = (10, 20, 50)$. *Numbers in the bracket are mean absolute errors of estimates*

| $n$ | True parameters | | MPS estimates | | MLEs | |
|---|---|---|---|---|---|---|
| | $\gamma_0$ | $\mu_0$ | $\check{\gamma}$ | $\check{\mu}$ | $\hat{\gamma}$ | $\hat{\mu}$ |
| 10 | −0.2 | 1 | −0.10(0.48) | 0.97(0.44) | −0.52(170) | 0.80(247) |
| | 0.2 | 1 | 0.41(0.55) | 0.92(0.48) | −0.18(7547) | 0.99(12347) |
| | 1 | 1 | 1.33(0.75) | 0.86(0.57) | 0.72(5304) | 1.20(17860) |
| | 1.2 | 1 | 1.58(0.81) | 0.83(0.59) | 0.92(3890) | 1.17(14582) |
| 20 | −0.2 | 1 | −0.13(0.26) | 0.98(0.28) | −0.45(0.27) | 0.97(0.28) |
| | 0.2 | 1 | 0.33(0.31) | 0.95(0.31) | −0.03(506) | 1.06(950) |
| | 1 | 1 | 1.21(0.46) | 0.91(0.37) | 0.87(70) | 1.08(300) |
| | 1.2 | 1 | 1.43(0.50) | 0.90(0.39) | 1.06(10) | 1.07(50) |
| 50 | −0.2 | 1 | −0.15(0.13) | 0.98(0.16) | −0.38(50) | 0.99(24) |
| | 0.2 | 1 | 0.27(0.18) | 0.96(0.21) | 0.07(50) | 1.02(50) |
| | 1 | 1 | 1.11(0.26) | 0.95(0.23) | 0.95(0.23) | 1.02(0.24) |
| | 1.2 | 1 | 1.32(0.28) | 0.94(0.24) | 1.14(0.25) | 1.03(0.25) |

We have also evaluated the empirical significance level of models having different

TABLE 3
*Failure rate of MPS estimation, maximum likelihood estimation and PWM estimation for the GEV Distribution. Tabulated values are the number of outliers per 100 simulated samples*

| $n$ | MPS estimation | | | | maximum likelihood estimation | | | | PWM estimation | | | |
|---|---|---|---|---|---|---|---|---|---|---|---|---|
| | $\gamma_0$ | | | | $\gamma_0$ | | | | $\gamma_0$ | | | |
| | -0.2 | 0.2 | 1 | 1.2 | -0.2 | 0.2 | 1 | 1.2 | -0.2 | 0.2 | 1 | 1.2 |
| 10 | 0.00 | 0.00 | 0.03 | 0.05 | 2.00 | 0.77 | 0.11 | 0.02 | 0.00 | 0.00 | 0.00 | 0.00 |
| 20 | 0.00 | 0.00 | 0.06 | 0.11 | 0.06 | 0.03 | 0.04 | 0.03 | 0.00 | 0.00 | 0.00 | 0.00 |
| 50 | 0.00 | 0.00 | 0.05 | 0.00 | 0.00 | 0.00 | 0.04 | 0.00 | 0.00 | 0.00 | 0.00 | 0.00 |

TABLE 4
*Failure rate of MPS estimation and maximum likelihood estimation for the GPD distribution. Tabulated values are the number of outliers per 100 simulated samples*

| $n$ | MPS estimation | | | | maximum likelihood estimation | | | |
|---|---|---|---|---|---|---|---|---|
| | $\gamma_0$ | | | | $\gamma_0$ | | | |
| | -0.2 | 0.2 | 1 | 1.2 | -0.2 | 0.2 | 1 | 1.2 |
| 10 | 0.00 | 0.00 | 0.00 | 0.00 | 0.05 | 2.51 | 3.61 | 2.94 |
| 20 | 0.00 | 0.00 | 0.00 | 0.00 | 0.00 | 0.19 | 0.06 | 0.01 |
| 50 | 0.00 | 0.01 | 0.00 | 0.00 | 0.01 | 0.01 | 0.00 | 0.00 |

TABLE 5
*Empirical sizes of Moran test statistics on the GEV from 10000 simulations of sample sizes $n = (10, 20, 50)$*

| $n$ | GEV Models | | | Empirical sizes | | |
|---|---|---|---|---|---|---|
| | $\gamma_0$ | $\mu_0$ | $\sigma_0$ | $\alpha = 0.10$ | $\alpha = 0.05$ | $\alpha = 0.01$ |
| 10 | $-0.2$ | 1 | 1 | 0.0618 | 0.0270 | 0.0034 |
| | 0.2 | 1 | 1 | 0.0580 | 0.0257 | 0.0037 |
| | 1 | 1 | 1 | 0.0592 | 0.0264 | 0.0062 |
| | 1.2 | 1 | 1 | 0.0619 | 0.0318 | 0.0115 |
| 20 | $-0.2$ | 1 | 1 | 0.0759 | 0.0337 | 0.0053 |
| | 0.2 | 1 | 1 | 0.0783 | 0.0373 | 0.0081 |
| | 1 | 1 | 1 | 0.0785 | 0.0360 | 0.0086 |
| | 1.2 | 1 | 1 | 0.0814 | 0.0377 | 0.0089 |
| 50 | $-0.2$ | 1 | 1 | 0.0848 | 0.0408 | 0.0077 |
| | 0.2 | 1 | 1 | 0.0906 | 0.0414 | 0.0074 |
| | 1 | 1 | 1 | 0.0890 | 0.0419 | 0.0101 |
| | 1.2 | 1 | 1 | 0.1000 | 0.0509 | 0.0143 |

TABLE 6
*Empirical sizes of Moran test statistics on the GPD from 10000 simulations of sample sizes $n = (10, 20, 50)$*

| $n$ | GPD Models | | Empirical sizes | | |
|---|---|---|---|---|---|
| | $\gamma_0$ | $\sigma_0$ | $\alpha = 0.10$ | $\alpha = 0.05$ | $\alpha = 0.01$ |
| 10 | -0.20 | 1 | 0.0745 | 0.0326 | 0.0044 |
| | 0.20 | 1 | 0.0699 | 0.0326 | 0.0044 |
| | 1.00 | 1 | 0.0718 | 0.0326 | 0.0044 |
| | 1.20 | 1 | 0.0734 | 0.0326 | 0.0044 |
| 20 | -0.20 | 1 | 0.0855 | 0.0394 | 0.0070 |
| | 0.20 | 1 | 0.0827 | 0.0374 | 0.0066 |
| | 1.00 | 1 | 0.0794 | 0.0375 | 0.0078 |
| | 1.20 | 1 | 0.0806 | 0.0377 | 0.0071 |
| 50 | -0.20 | 1 | 0.0972 | 0.0474 | 0.0093 |
| | 0.20 | 1 | 0.0932 | 0.0452 | 0.0081 |
| | 1.00 | 1 | 0.0960 | 0.0471 | 0.0091 |
| | 1.20 | 1 | 0.0940 | 0.0477 | 0.0087 |

TABLE 7

*GEV parameter estimation by MPS in* 10000 *simulations of sample sizes* $n = 10, 20, 50$. *Data are generated from exponential distributions with* $\lambda = 0.1, 0.5, 1.0, 5.0$. *Figures shown are* 25%, 50% *and* 75% *quantiles of the estimated parameter.*

|   |   | Parameter quantile estimates of GEV ||||||||
|---|---|---|---|---|---|---|---|---|---|
| $n$ | $\lambda$ | 25%$\breve{\gamma}$ | 50%$\breve{\gamma}$ | 75%$\breve{\gamma}$ | 25%$\breve{\mu}$ | 50%$\breve{\mu}$ | 75%$\breve{\mu}$ | 25%$\breve{\sigma}$ | 50%$\breve{\sigma}$ | 75%$\breve{\sigma}$ |
| 10 | 0.1 | −0.33 | −0.07 | 0.22 | 31.51 | 33.82 | 36.27 | 8.33 | 10.56 | 12.88 |
|    | 0.5 | −0.33 | −0.07 | 0.21 | 6.3   | 6.76  | 7.25  | 1.67 | 2.11  | 2.58  |
|    | 1.0 | −0.33 | −0.07 | 0.22 | 3.15  | 3.38  | 3.63  | 0.83 | 1.06  | 1.29  |
|    | 5.0 | −0.32 | −0.04 | 0.27 | 0.62  | 0.68  | 0.73  | 0.17 | 0.22  | 0.28  |
| 20 | 0.1 | −0.18 | −0.05 | 0.11 | 32.28 | 33.93 | 35.65 | 8.97 | 10.39 | 11.79 |
|    | 0.5 | −0.19 | −0.05 | 0.11 | 6.46  | 6.79  | 7.13  | 1.79 | 2.08  | 2.36  |
|    | 1.0 | −0.19 | −0.05 | 0.11 | 3.23  | 3.39  | 3.57  | 0.9  | 1.04  | 1.18  |
|    | 5.0 | −0.19 | −0.05 | 0.11 | 0.65  | 0.68  | 0.71  | 0.18 | 0.21  | 0.24  |
| 50 | 0.1 | −0.11 | −0.03 | 0.05 | 33.01 | 34.06 | 35.14 | 9.33 | 10.13 | 10.99 |
|    | 0.5 | −0.11 | −0.03 | 0.05 | 6.6   | 6.81  | 7.03  | 1.87 | 2.03  | 2.2   |
|    | 1.0 | −0.11 | −0.03 | 0.05 | 3.3   | 3.41  | 3.51  | 0.93 | 1.01  | 1.1   |
|    | 5.0 | −0.11 | −0.03 | 0.05 | 0.66  | 0.68  | 0.70  | 0.19 | 0.2   | 0.22  |

scale and location parameters. The results did not differ significantly and thus were not reported here. It seems that the performance of Moran's statistic was affected by the sample size rather than the underlying models.

In application, it is common to take cluster maxima in the model fitting of the GEV. Having shown that the MPS gives stable estimations on data generated from known models, in the following, fitting the maximum observations in clusters of size 30 was performed. This experiment mimics the situation that the original data are daily observations with GEV fitted to the monthly maxima. The aim of this experiment is to evaluate the stability in the estimation of cluster maxima.

Data $x_{n,m}$ were simulated from the exponential distribution

$$F(x) = 1 - e^{-\lambda x} \quad x > 0$$

with $\lambda = 0.1, 0.5, 1.0, 5.0$ where $n$ was the sample size of maxima and $m = 30$ the size of a cluster. From each cluster, the maximum, $\max(x_{n,1}, \ldots, x_{n,30})$, was taken and the GEV distributions was fitted to the data by MPS method.

Table 7 shows the estimated parameter quantiles. In the GEV fitting, the tail estimates fall in a narrow range in the four cases $\lambda = 0.1, 0.5, 1.0, 5.0$. Note that the medians for $\breve{\sigma}$ are proportional to the value of $\lambda^{-1}$. This feature remains stable across all sample sizes. A similar pattern is also observed for the medians of $\breve{\mu}$. This again shows that estimation using MPS is stable and reliable.

## 5. Real examples

Some real data sets were studied in the literature (Castillo et al., [1]) using the maximum likelihood method. To illustrate the advantages of the MPS approach, in this paper, four examples were studied, namely, the age data, the wave data, the wind data and the flood data. The above four data sets are obtainable in Castillo et al. [1]. The first example is the oldest age of men at death in Sweden. The annual oldest ages at death in Sweden from 1905 to 1958 were recorded. The age data may be used to predict oldest ages at death in the future. The wave data set contains the yearly maximum heights, in feet. The data could be used in the design of a breakwater. Then, in the wind data, the yearly maximum wind speed in miles per hour is considered. A wind speed design for structural building purposes could be determined from this data set. The last example is the flood data which consists of the yearly maximum flow discharge, in cubic meters. The data may help in designing a flood protection device.

TABLE 8
*Estimated GPD parameters by MPS in four examples*

| Data | Threshold | $\breve{\gamma}$ | $\breve{\sigma}$ | $M(\breve{\boldsymbol{\theta}})$ |
|---|---|---|---|---|
| Age | 104.01 | 1.06 | 2.79 | 43.01 |
| Wave | 17.36 | 0.01 | 7.00 | 45.81 |
| Flood | 45.04 | −0.03 | 9.62 | 60.53 |
| Wind | 36.82 | −0.88 | 5.31 | 46.79 |

TABLE 9
*Estimated GPD parameters by maximum likelihood method in four examples*

| Data | Threshold | $\hat{\gamma}$ | $\hat{\sigma}$ | Log-likelihood |
|---|---|---|---|---|
| Age | 104.01 | 1.38 | 3.45 | −9.23 |
| Wave | 17.36 | 0.27 | 7.98 | −39.33 |
| Flood | 45.04 | 0.20 | 10.87 | −57.34 |
| Wind | 36.82 | −0.48 | 6.52 | −47.01 |

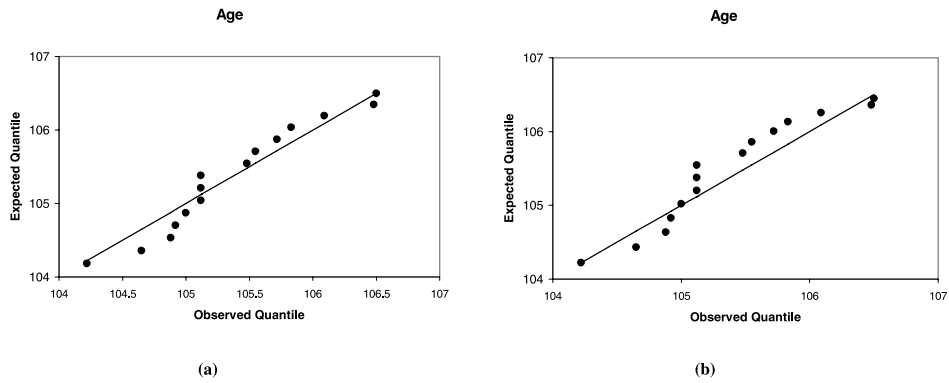

Fig 1. *Quantile plots of the age data fitted with the GPD using the MPS* (a) *and the MLE* (b).

In this section, we focus on the GPD with the maximum likelihood method and the MPS method. The GPD was fitted to the excess over a threshold. The thresholds were taken from [1]. Fitted parameters are shown in Tables 8 and 9. Note that $\breve{\gamma}$ and $\hat{\gamma}$ are greater than 1 for the age data. They are less than 1 for the wave, flood and wind data sets.

### 5.1. The GPD model for age data

Recall from Theorem 3.1 that the MLE does not exist for $\gamma > 1$. When the GPD is fitted to the age data, maximization of the GPD log-likelihood leads to the estimate $(\hat{\gamma}, \hat{\sigma}) = (1.38, 3.45)$ for which the log-likelihood is $-9.23$. The corresponding values using MPS are $(\breve{\gamma}, \breve{\sigma}) = (1.06, 2.79)$ and $M(\breve{\boldsymbol{\theta}}) = 43.01$. Fig. 1 shows the quantile plot of the two models fitted by the MPS and the maximum likelihood method respectively. In each plot, the expected quantile is calculated by

$$\text{GPD} : Q_{\text{GPD}} = \left[1 - \left(1 - \frac{i}{n+1}\right)^{\gamma}\right] \frac{\sigma}{\gamma} \quad (i = 1, 2, \ldots, n).$$

where $\mathbf{\Theta}_{\text{GPD}} = (\gamma, \sigma)^T$ are estimated parameters either by the MPS or by the maximum likelihood method.

The MPS seems to perform better than the MLE. Empirical upper quantiles in the MPS are closer to that of a straight line. This suggests that the MPS is a better method in this case.

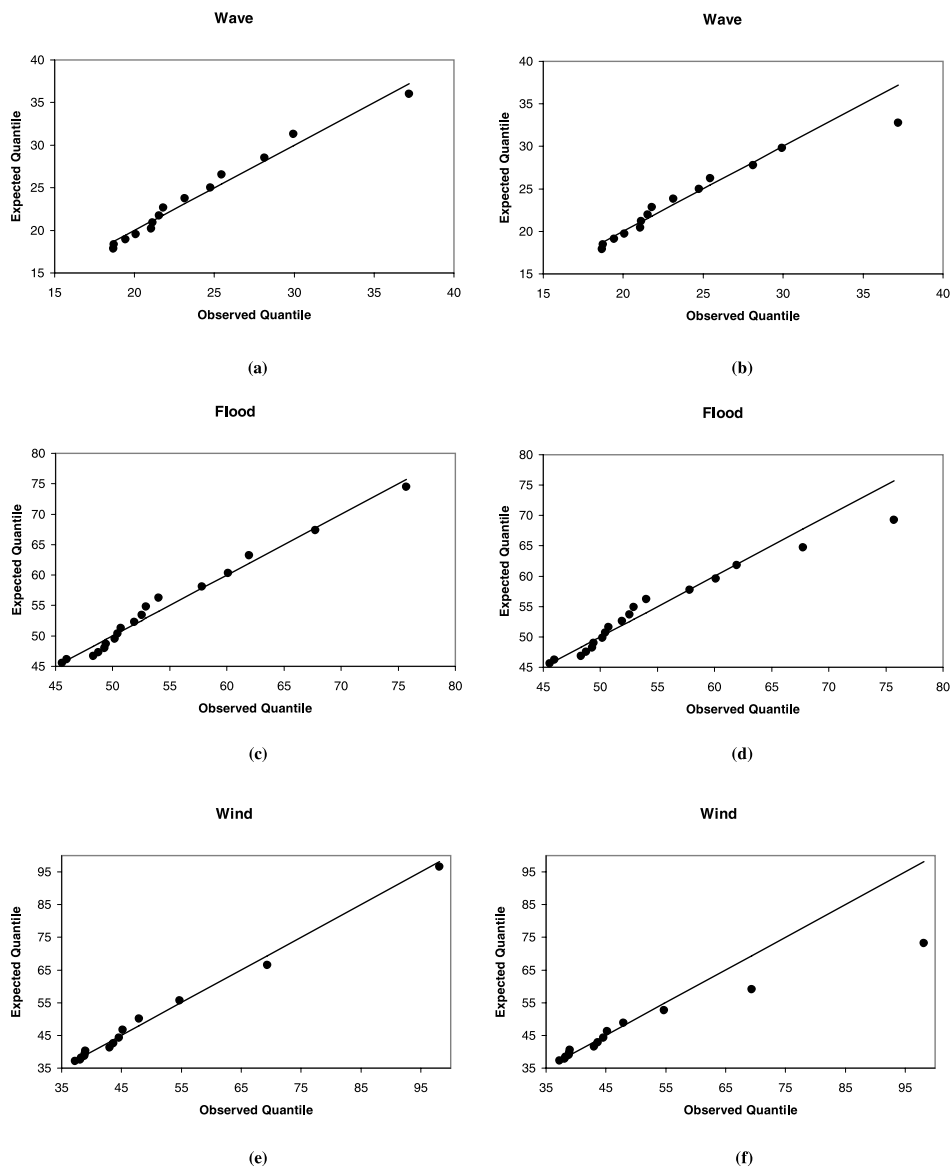

FIG 2. *Quantile plots of the wave data ((*a*) and (*b*)), the flood data ((*c*) and (*d*)) and the wind data ((*e*) and (*f*)) fitted with the GPD. The expected quantiles of (*a*), (*c*) and (*e*) were based on the MPS. The expected quantiles of (*b*), (*d*) and (*f*) were based on the MLE.*

### 5.2. GPD model for the wave data, flood data and wind data

The GPD was also considered for the wave data, the flood data and the wind data. Thresholds for the GPDs were taken as in Castillo et al. [1]. The quantile pots for the MPS are reported in Fig. 2(a), 2(c) and 2(e) and those for the MLE are reported in Fig. 2(b), 2(d) and 2(f). With reference to Fig. 2(a), 2(c) and 2(e), it can be seen that empirical quantiles based on the MPS keep close to the fitted model's. However, in Fig. 2(b), 2(d) and 2(f), plots of the upper quantiles based on the MLE seem to deviate more from a straight line. This suggests that the MPS gives a better fit to the data.

## 6. Conclusion and discussion

In extreme value analysis, one technical problem is the lack of data owing to the fact that only extreme observations are used for model fitting. Subject to this constraint, a method that is able to give stable estimates is highly desirable. Juarez and Schucany [9] have demonstrated the instability of the influence curve of the MLE at small sample sizes. This is in agreement with the presented simulation results. In contrast, the MPS works satisfactorily. Not only does the MPS yield closer estimates from data generated from a known parameter set, it also keeps performing stably for data maxima taken from clusters. It also works well under $\gamma \geq 1$ whereas the MLE does not. In addition to MPS's simple formulation and execution, its by-product, the Moran's statistic, is shown to perform well in checking the goodness of fit. The MPS could potentially be one of the best methods in fitting extreme value distributions. On the other hand, it has been shown in [2] that the MPS is a function of sufficient statistics. Extension to multivariate problems using MPS is also going to be explored.

### Acknowledgements

W. K. Li would like to thank the organiser of the C. Z. Wei Memorial Conference for the invitation to participate in the conference. Partial support by the Area of Excellence Scheme under the University Grants Committee of the Hong Kong Special Administration Region, China (Project AoE/P-04/2004) is also gratefully acknowledged. The authors thank a referee for comments that led to improvement of the paper.